\documentclass{amsart}

\usepackage{amsmath,amssymb}
\usepackage{amsthm}
\usepackage{amsrefs}
\usepackage{indentfirst}
\usepackage{mathrsfs}
\usepackage{graphicx}
\usepackage{hyperref}
\usepackage{xcolor}
\usepackage{fourier}
\usepackage[margin=1.5in]{geometry}
\usepackage[norelsize,boxed,linesnumbered,vlined,algo2e]{algorithm2e}

\theoremstyle{plain}

\theoremstyle{definition}

\theoremstyle{remark}

\DeclareMathOperator{\diag}{diag}

\newcommand{\ud}{\,\mathrm{d}}

\newcommand{\ZZ}{\mathbb{Z}}
\newcommand{\TT}{\mathbb{T}}

\newcommand{\wh}[1]{\widehat{#1}}

\IfFileExists{mathabx.sty}%
  {\DeclareFontFamily{U}{mathx}{\hyphenchar\font45}%
   \DeclareFontShape{U}{mathx}{m}{n}{<->mathx10}{}%
   \DeclareSymbolFont{mathx}{U}{mathx}{m}{n}%
   \DeclareFontSubstitution{U}{mathx}{m}{n}%
   \DeclareMathAccent{\widebar}{0}{mathx}{"73}%
}{%
  \PackageWarning{mathabx}{%
    Package mathabx not available, therefore\MessageBreak substituting
    widebar with overline\MessageBreak }%
  \newcommand{\widebar}[1]{\overline{#1}}%
}

\newcommand{\mc}[1]{\mathcal{#1}}

\newcommand{\veps}{\varepsilon}

\newcommand{\abs}[1]{\lvert#1\rvert}
\newcommand{\Abs}[1]{\left\lvert#1\right\rvert}

\newcommand{\norm}[1]{\lVert#1\rVert}

\newcommand{\average}[1]{\left\langle#1\right\rangle}

\title{Sparsifying preconditioner for soliton calculations}
\author{Jianfeng Lu} \address{Department of Mathematics, Department of
  Physics, Department of Chemistry, and Fitzpatrick Institute of
  Photonics, Duke University, Box 90320, Durham NC 27708, USA}
\email{jianfeng@math.duke.edu}

\author{Lexing Ying}
\address{Department of Mathematics and Institute of Computational and Mathematical Engineering, Stanford University, 450 Serra Mall, Bldg 380, Stanford CA 94305, USA}
\email{lexing@stanford.edu}

\date{\today} \thanks{The work of J.L. is supported in part by the
  Alfred P.~Sloan Foundation and the National Science Foundation under
  grant DMS-1312659 and DMS-1454939. The work of L.Y. is partially
  supported by the National Science Foundation under grant DMS-0846501
  and the U.S. Department of Energy's Advanced Scientific Computing
  Research program under grant
  DE-FC02-13ER26134/DE-SC0009409. J.L. would like to thank helpful
  discussions with Jeremy Marzuola and Mikael Rechtsman.}

\begin{document}

\begin{abstract}
  We develop a robust and efficient method for soliton calculations
  for nonlinear Schr\"odinger equations. The method is based on the
  recently developed sparsifying preconditioner combined with Newton's
  iterative method. The performance of the method is demonstrated by
  numerical examples of gap solitons in the context of nonlinear
  optics.
\end{abstract}

\maketitle

\section{Introduction}

\subsection{Background} 

In this work, we consider a type of nonlinear equations associated
with nonlinear eigenvalue problems: Given $\lambda$, we solve for $u$
such that
\begin{equation}\label{eq:noneig}
  - \Delta u + Vu + N(x, u) = \lambda u
\end{equation}
with periodic boundary condition on a box domain $\TT^d := [0, 1]^d$,
where $d$ is the dimension.  Here $V$ is some background potential and
$N$ is the nonlinear term.  This is a natural generalization of linear
eigenvalue problem, for instance
\begin{equation}\label{eq:lineig}
  -\Delta u + V u = \lambda u
\end{equation}
with periodic boundary condition. Note that while for the linear
eigenvalue problem \eqref{eq:lineig}, a nontrivial solution $u$ exists
only when $\lambda$ is an eigenvalue and the normalization of $u$ is
not fixed (the eigenvalue remains the same and the eigenfunction is
simply scaled), the solution to the nonlinear problem
\eqref{eq:noneig} will have a fixed normalization due to the presence
of the nonlinear term. In fact, for most cases, we are interested in
the relation between $\lambda$ (which can be seen as a Lagrange
multiplier) and the normalization $\norm{u}_2$. 

This type of nonlinear eigenvalue problems arises in several
applications, for example in nonlinear photonic crystals (see e.g.,
\cites{Efremidis:03, SoljacicJoannopoulos:04, SlusherEggleton:04,
  DudleyTaylor:09} and the numerical examples in
Section~\ref{sec:numerics}) and in the study of solitons in nonlinear
wave equations, for instance the Gross-Pitaevskii equation. It is also
related to electronic structure theory, for which a system of
nonlinear eigenvalue problem arises as the Kohn-Sham equations. While
we will focus on applications in nonlinear Schr\"odinger equations for
solitons in photonic crystals in this work, the algorithm applies to
other context as well.

While this paper considers numerical methods to solve the equation
\eqref{eq:noneig}, let us mention that the analytical understanding of
such equations is also a very active research area in applied
mathematics, see e.g., a recent review \cite{Weinstein:15}.

\subsection{Previous approaches}

To get the nonlinear eigenvalue $\lambda$ and the corresponding state
$u$, one popular approach in the engineering literature is to use an
indirect strategy based on time evolution. To get an estimate of the
eigenvalue, the time-dependent equation
\begin{equation}
  i u_t = \Delta u - V(u) - N(x, u)
\end{equation}
is evolved for a given time period with a random initial data; a
Fourier transform in time is then used to extract the frequencies with
significant amplification at the final time compared to the initial
data. While the time evolution method is very useful to extract
resonant frequencies for linear problems, this procedure becomes quite
tricky if the nonlinear effects become significant (see e.g.,
\cite{Tran:95}). In particular, the lack of the superposition
principle makes it difficult to interpret the amplification in the
Fourier space.

Thus, in the presence of nonlinearity, it is preferred to solve the
nonlinear eigenvalue problem directly (without resorting to
time-evolution) based on iterative methods. One such approach is the
self-consistent iteration scheme based on formulating the problem as a
linear eigenvalue problem at each iteration, e.g., 
\begin{equation}
  -\Delta u^{(n+1)} + V u^{(n+1)} + \frac{N(x, u^{(n)})}{u^{(n)}} u^{(n+1)} = \lambda^{(n+1)} u^{(n+1)}
\end{equation}
by choosing $\lambda^{(n+1)}$ the closest eigenvalue to the desired
$\lambda$. The difficulty of this type of method is that we do not
have a priori knowledge on the normalization of $u^{(n+1)}$ (recall
that the normalization is arbitrary at the linear level, but has to be
fixed for the nonlinear equation), so that some heuristics or
parameter tuning is needed.

Another more popular iterative scheme was proposed originally by
Petviashvili \cite{Petviashvili:76} for equation of the type
\begin{equation}
  -\Delta u - u^3 = \lambda u 
\end{equation}
for $\lambda > 0$, so that the $N(x, u)$ is the focusing cubic
nonlinearity $-u^3$ and there is no potential $V$. The iteration
scheme reads
\begin{equation}
  \wh{u^{(n+1)}}(k) = M_n^{\gamma} \frac{\mc{F}\bigl((u^{(n)})^3\bigr)(k)}{\lambda + 4\pi^2 \abs{k}^2}, 
\end{equation}
where $M_n$ is a stabilizing factor given by 
\begin{equation}
  M_n = \dfrac{\displaystyle \sum_k (\lambda + 4\pi^2 \abs{k}^2)\abs{\wh{u^{(n)}(k)}}^2}{\displaystyle \sum_k \wh{u^{(n)}(k)} \mc{F}\bigl( (u^{(n)})^3 \bigr)(k)},
\end{equation}
and $\gamma$ is a parameter to be chosen so that the iteration
converges. Here $\wh{u}$ and $\mc{F}[u]$ both denote the Fourier
transform of $u$: For each $k \in \ZZ^d$,
\begin{equation}
  \wh{u}(k) = \mc{F}[u](k) := \int_{\TT^d} e^{-2\pi i (k \cdot x)} u(x) \ud x. 
\end{equation}
The convergence of the Pitviashvili scheme has been analyzed in
\cite{PelinovskyStepanyants:04}. The iteration scheme was also
generalized to deal with more general nonlinearities and the case with
potential terms, see for example \cite{AblowitzMusslimani:05,
  LakobaYang:07} and references therein.  All methods of this type
reformulate the nonlinear problem \eqref{eq:noneig} using the Fourier
space representation, in which the positively shifted Laplacian
operator becomes diagonal, and hence can be explicitly inverted.

\subsection{Our contribution}

In this work, to solve the nonlinear equation \eqref{eq:noneig}, we
will use the standard Newton iteration method. Given $\lambda$, to
find the eigenfunction $u$ (and hence its normalization), the Newton
iterative scheme is given by
\begin{equation}
  u^{(n+1)} = 
  u^{(n)} - \bigl(-\Delta + L_{u^{(n)}} - \lambda \bigr)^{-1} 
  \bigl(
    - \Delta u^{(n)} + V u^{(n)} + N(x, u^{(n)}) - \lambda u^{(n)} \bigr).
\end{equation}
Here $L_u$ is the linearization of the potential terms at $u$: $L_u =
V + \frac{\delta N}{\delta u}$. The key step for the above iteration
scheme is the solution to the linear system
\begin{equation}\label{eq:linear}
  \bigl(-\Delta + L_{u} - \lambda\bigr) v = r
\end{equation}
for the current iterate $u$ and remainder denoted as $r$.

We will consider the Fourier pseudospectral method to discretize
\eqref{eq:noneig}, which is very popular for the nonlinear wave
equations with periodic boundary conditions. The pseudospectral method
typically requires minimal degree of freedoms for a given accuracy
among standard discretizations and is also easy to implement, and
hence widely used in physics and engineering literature (see e.g.,
\cite{Plihal:94} in the context of photonics). As we have discussed
already, most previous iterative methods are built on the Fourier
pseudospectral method.

The drawback of the pseudo-spectral method however is that the
discrete matrices resulting from discretizing $-\Delta + L_u -
\lambda$ are dense, so that it is usually very expensive to directly
solve equations like \eqref{eq:linear}. Therefore, some iterative
methods have to be used, in which a shifted Laplacian is often used as
a preconditioner. In fact, the iterative methods introduced by
Petviashvili \cite{Petviashvili:76} and the various generalizations
can be understood as using the shifted Laplacian in a way that the
outer and inner iterations are combined.

The main motivation of this work is that, rather than finding an
alternative iterative scheme which usually is less effective than the
Newton iteration, we will instead use an efficient way to solve the
dense discrete system resulting from \eqref{eq:linear}. The key
component is the recently developed sparsifying preconditioner
\cite{Ying:spspd}, which transforms the discrete dense systems
numerically to a sparse one. Sparse direct solvers can then be
combined with standard iterative schemes to provide an efficient way
to solve the linear system \eqref{eq:linear}, and hence the nonlinear
equation.

\section{Algorithm description}

\subsection{Sparsifying preconditioner}\label{sec:sparse}

Let us focus on numerical solution of \eqref{eq:linear}, given the
current iterate $u$ and hence $L_u$. We assume the nonlinearity term
$N(x, u)$ only depends locally on $u$, and hence the derivative
$\frac{\delta N}{\delta u}$ with respect to $u$ gives a potential-type 
term in \eqref{eq:linear}.

We assume that the computation domain is the periodic unit box
$[0,1]^d$ and discretize the problem using the Fourier pseudospectral
method. Let us index the set of all grid points as
\begin{equation}
  J = \bigl\{(j_1, \ldots, j_d) \mid 0 \leq j_1, \ldots, j_d < n \bigr\}. 
\end{equation}
For each $j \in J$, in the pseudospectral method, the grid point is given by 
\begin{equation}
  x_j = jh
\end{equation}
with $h = \frac{1}{n}$ being the mesh size in each dimension. The
function $v$ (and similarly $u$ and $r$, etc.) is hence discretized as
vectors with
\begin{equation}\label{eq:v}
  v_j = v(jh). 
\end{equation}
The corresponding Fourier grid is then given by 
\begin{equation}
  K = \bigl\{ (k_1, \ldots, k_d) \mid -n/2 \leq k_1, \ldots, k_d < n/2 \bigr\}. 
\end{equation}
In the pseudospectral method, the Laplacian operator is discretized as
(with slight abuse of notation, we still denote by $-\Delta$ the
discretized operator)
\begin{equation}\label{eq:Delta}
  - \Delta = F^{-1} \diag(4 \pi^2 \abs{k}^2)_{k \in K} F, 
\end{equation}
where $F$ and $F^{-1}$ are the discrete Fourier and inverse Fourier
transforms
\begin{align}
  (F f)_k & = \frac{1}{n^d} \sum_{j \in J} e^{-2\pi i (k \cdot j) / n} f_j, \qquad k \in K; \\
  (F^{-1} g)_j & = \sum_{k\in K} e^{2\pi i (j \cdot k) / n} g_k, \qquad j \in J. 
\end{align}
So that $-\Delta$ is diagonal in the Fourier space. The discretized
$L_u$ (again, we use the same notation after discretization), on the
other hand, is a diagonal matrix in the physical space,
\begin{equation}\label{eq:Lu}
  L_u = \diag\bigl(L_u(jh)\bigr)_{j \in J}. 
\end{equation}
Following the notations introduced in in \eqref{eq:v},
\eqref{eq:Delta}, and \eqref{eq:Lu}, the discretized problem of
\eqref{eq:linear} takes the same form
\begin{equation}\label{eq:Disc}
  \bigl(-\Delta + L_{u} - \lambda\bigr) v = r.
\end{equation}
We will briefly recall the key idea of the sparsifying preconditioner
for solving this type of equations. More details can be found in
\cite{Ying:spspd} (see also \cite{Ying:spspc}).

Denote $l = \average{L_u}$ the spatial average of $L_u$. We assume
without loss of generality that $-\Delta + l - \lambda$ is invertible
on $\TT^d$ with periodic boundary condition, otherwise, we will use a
slight perturbation of $\lambda$ instead and put the difference into
$L_u$. This allows us to rewrite \eqref{eq:linear} trivially as
\begin{equation}
  \bigl(-\Delta + (l-\lambda)  + (L_{u}-l) \bigr) v = r
\end{equation}
The inverse $G$ of the constant part $ \bigl(-\Delta + (l-\lambda)
\bigr)$ is given by
\[
G = F^{-1} \diag\bigl(\frac{1}{4 \pi^2
  \abs{k}^2+l-\lambda}\bigr)_{k \in K} F.
\]
Applying $G$ to both sides of \eqref{eq:linear} gives an equivalent
linear system
\begin{equation}\label{eq:GLu}
  \bigl(I + G (L_u - l) \bigr) v = G r. 
\end{equation}
Note that in the pseudo-spectral discretization, since $G$ can be
easily evaluated by fast Fourier transform and $L_u - l$ is a diagonal
matrix, applying the operators on both sides of \eqref{eq:GLu} can be
carried out efficiently. The main idea of the sparsifying
preconditioner is to introduce a special sparse matrix $Q$ and
multiply $Q$ on the both hand sides of \eqref{eq:GLu}:
\begin{equation}\label{eq:QG}
  (Q + QG (L_u - l) ) v = Q (Gr).
\end{equation}
The operator $Q$ is required to satisfy two conditions:
\begin{itemize}
\item For each point $j$, the row $Q(j,:)$ is supported in a local
  neighborhood $\mu(j)$ (to be defined below);
\item For each point $j$, $(QG)(j,\mu(j)^c)=Q(j,\mu(j))G(\mu(j),\mu(j)^c) \approx 0$.
\end{itemize}
Let $S =\{(j,\mu(j)), j\in J\}$ be the support of the sparse matrix
$Q$. The second condition ensures that the product $QG$ is also
essentially supported on $S$. The sum $Q+QG (L_u-l)$ is essentially
supported in $S$ as well, since $(L_u-l)$ is a diagonal matrix. Define
$P$ to be the restriction of $Q+QG (L_u-l)$ to $S$ by thresholding
other values to zero, i.e.
\begin{equation}\label{eq:P}
  P_{ij}=
  \begin{cases}
    \bigl(Q+QG (L_u-l)\bigr)_{ij}, & (i,j)\in S,\\
    0, &  (i,j)\not\in S.
  \end{cases}
\end{equation}
Knowing that $P$ is a close approximation of $Q+QG (L_u-l)$, one
arrives at the approximate equation
\begin{equation}
  P v \approx Q(Gr). 
\end{equation}
The sparsifying preconditioner computes an approximate solution
$\tilde{v}$ by solving
\[
P \tilde{v} = Q(Gr).
\]
Since $P$ is sparse, the above equation can be solved by sparse direct
solvers such as the nested dissection algorithm \cite{George:73}. The
nested dissection algorithm constructs an efficient sparse
factorization of $P^{-1}$ by exploiting the sparse pattern of
$P$. Once the factorization is ready, applying $P^{-1}$ to a given
factor can be done efficiently. The solution $\tilde{v}=P^{-1} QG r$
can be used as a preconditioner for the standard iterative algorithms
such as GMRES \cite{SaadSchultz:86} for the solution of
\eqref{eq:Disc}.

As detailed in \cite{Ying:spspd}, the neighborhood $\mu(j)$ in the
first condition is chosen to be the stencil support of a standard
spectral element method so that one can fully exploit the efficiency
of the nested dissection method used for computing the approximate
inverse for \eqref{eq:QG}.

\subsection{Newton iteration}

With the solution of \eqref{eq:linear} as discussed above, we solve
the nonlinear equation \eqref{eq:noneig} using the standard Newton
iteration. 
We summarize here the overall algorithm as
\begin{algorithm2e}[h]
\label{alg}
Initialize: Set $u^{(0)}$ and $n = 0$; 

\While{not converged}{
  Calculate $r^{(n)} = - \Delta u^{(n)} + V u^{(n)} + N\bigl(x, u^{(n)}\bigr) - \lambda u^{(n)}$; 

  Construct the sparsifying preconditioner $Q$ following the
  discussion in \S\ref{sec:sparse}; 

  Construct and form a nested dissection factorization of the sparse
  matrix $P$ in \eqref{eq:P};
  
  Solve $\bigl(-\Delta + (l-\lambda) + (L_{u}-l) \bigr) v^{(n)} =
  r^{(n)}$ for $v^{(n)}$ with $P^{-1}QG$ as the
  preconditioner;

  $u^{(n+1)} = u^{(n)} - v^{(n)}$; 

  Set $n = n + 1$. 
}
\end{algorithm2e}

\section{Numerical examples}\label{sec:numerics}

To test the Algorithm~\ref{alg}, we consider a few examples from the
optics, drawn from the paper \cite{Efremidis:03}.  Up to a rescaling
of the domain and various parameters and variables, these nonlinear
eigenvalue problem are of the type \eqref{eq:noneig}. The algorithms
presented above are implemented in MATLAB and the results reported in
this section are obtained on a Linux computer with a 2.6GHz CPU.

One example is Kerr nonlinearity, the equation is given by (using the
notations of \cite{Efremidis:03})
\begin{equation}\label{eq:kerr}
  -\frac{1}{2} \Delta u + \frac{V_0}{2} \bigl[\sin^2(\pi x) + \sin^2(\pi y)\bigr] u 
  - \sigma \abs{u}^2 u = \lambda u, 
\end{equation}
where $\sigma = \pm 1$ corresponds to focusing and defocusing cases.
The other example we consider is the case of saturable nonlinearity:
\begin{equation}\label{eq:satu}
  - \frac{1}{2} \Delta u + \frac{ V_0 u }{1 + A^2 \cos^2(\pi x)
    \cos^2(\pi y) + \abs{u}^2 } = \lambda u.
\end{equation}
In both examples, the computational domain is taken to be sufficiently
large in order to make sure that the calculation is not suffered from
finite size artifacts.

The particular interesting cases are when the eigenvalue (denoted by
$\lambda$ in above equations) lies in the band gaps of the linear
operator (the Bloch-Floquet theory applies thanks to the periodic
structure). In the literature, the solutions are referred as gap
solitons. From a mathematical point of view, this suggests that the
operator on the left hand side of \eqref{eq:linear} is indefinite,
which contains both positive and negative eigenvalues. Thus, the
negatively shifted Laplacian will not be a good preconditioner in this
case, which makes Pitviashvili type iterative methods very difficult
to converge. In comparison, the Newton iteration captures the correct
spectral behavior of the linearized operator, and is hence much more
robust.

\begin{figure}[ht!]
  \begin{tabular}{cc}
    \includegraphics[height=1.5in]{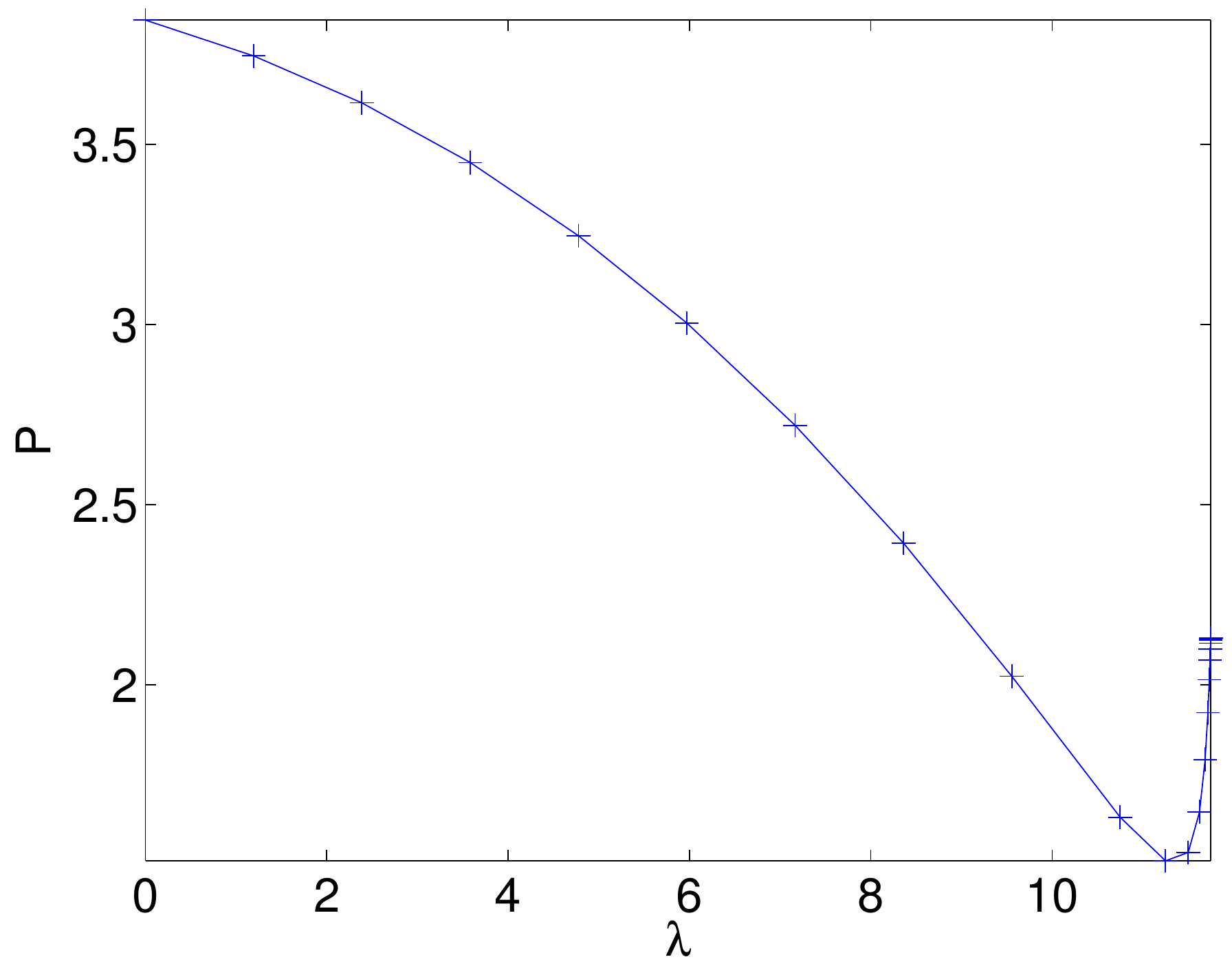} & \includegraphics[height=1.5in]{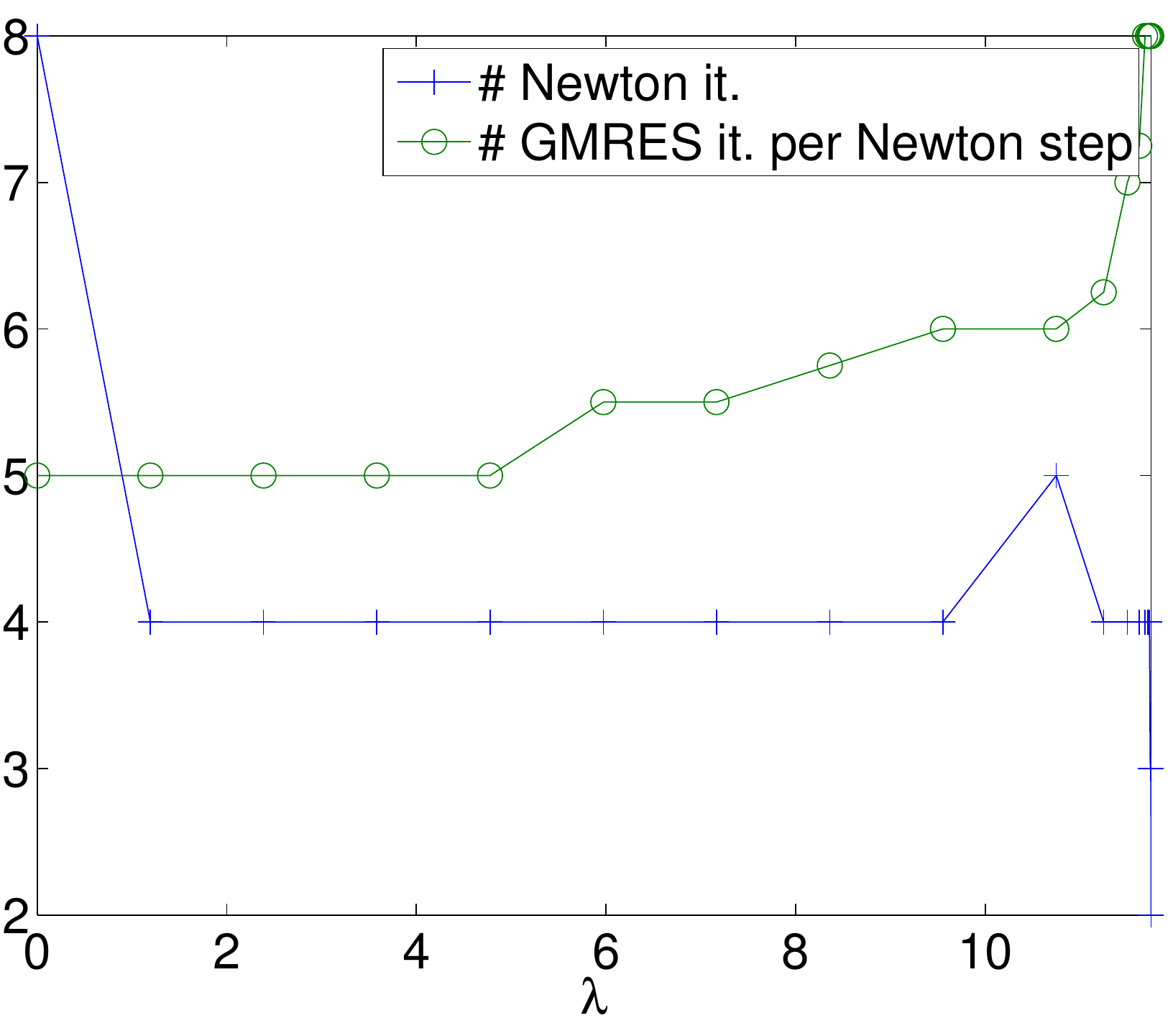}\\
    (a) & (b)\\
    \includegraphics[height=2in]{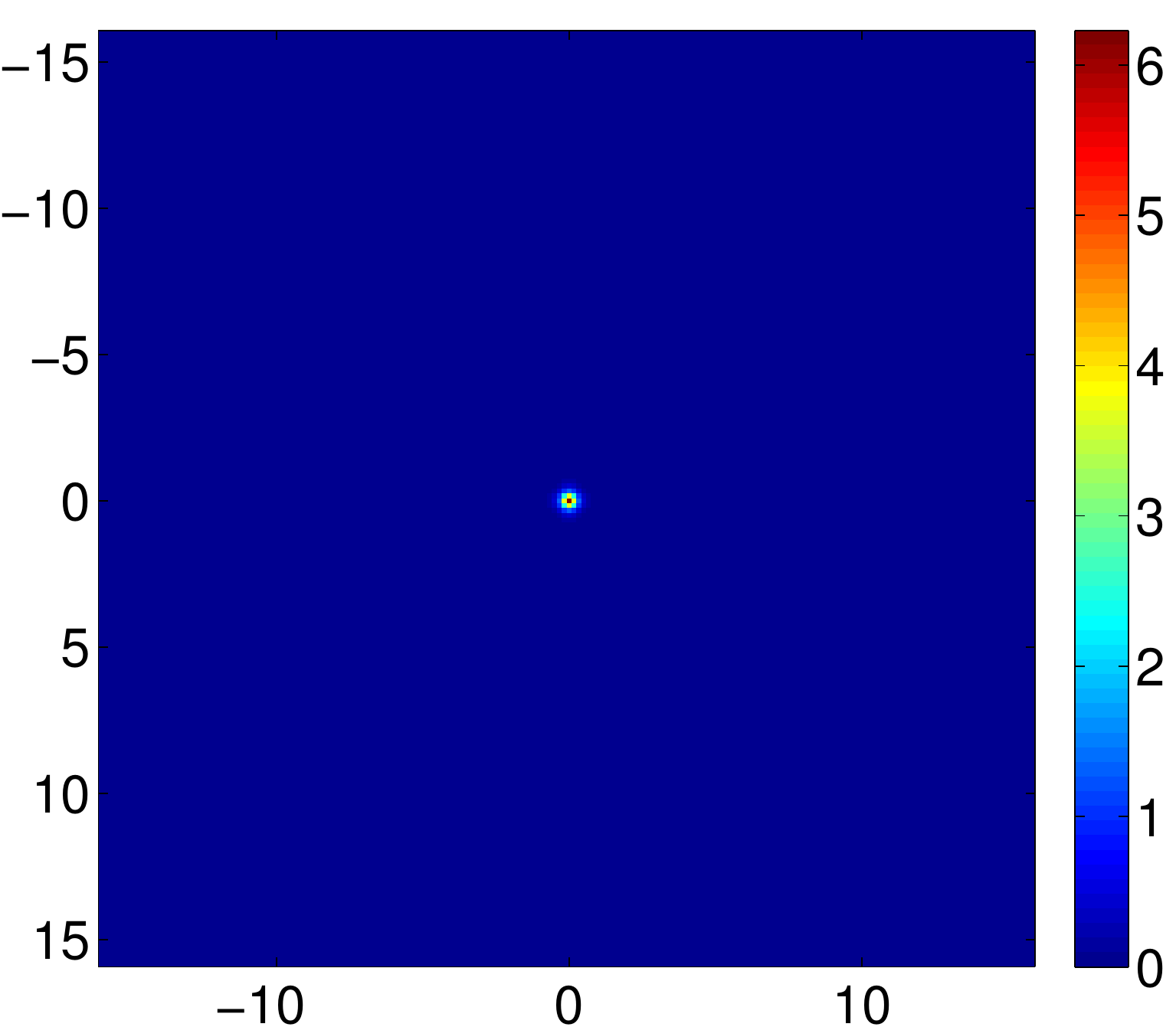} & \includegraphics[height=2in]{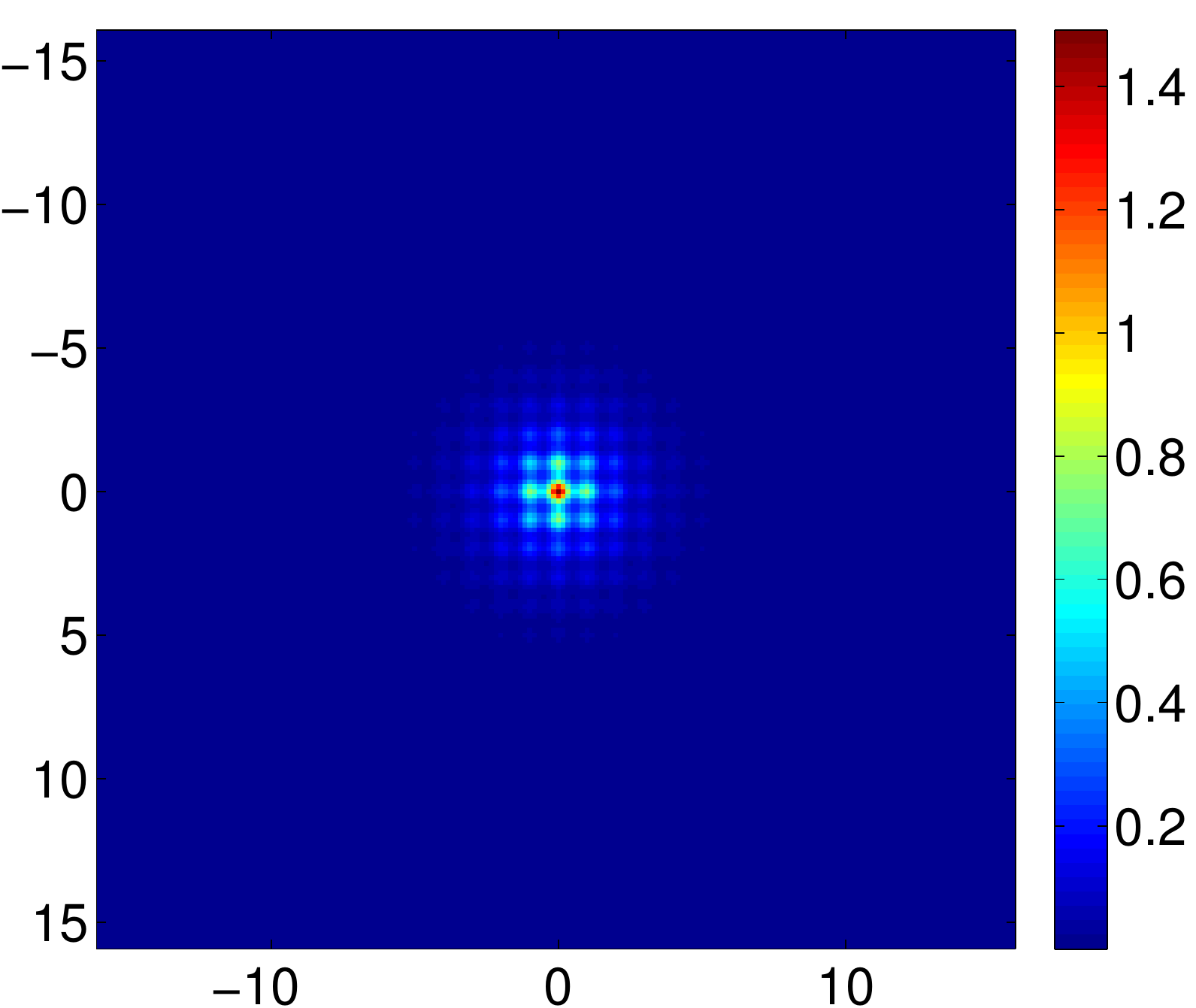}\\
    (c) & (d)
  \end{tabular}
  \caption{Kerr nonlinearity focusing case $V_0=28.8$. (a) the
    $\lambda-P$ curve. (b) The number of Newton iterations for each
    $\lambda$ value and the averaged number of preconditioned GMRES
    iteration for each Newton step.  (c) the field profile at
    $\lambda=0$. (d) the field profile at $\lambda=11.7498$.  }
  \label{fig:t4}
\end{figure}

Figure \ref{fig:t4} summarizes the results of \eqref{eq:kerr} in the
focusing case ($\sigma = 1$) with $V_0=28.8$. The computation is
carried out in the periodic domain $[-16,16]^2$ with $6$ points per
unit length. Figure \ref{fig:t4}(a) plots the relationship between the
eigenvalue $\lambda$ of the soliton and the power $P$ of the solution
defined by $\int |u(x)|^2 dx$. This $\lambda-P$ curve is obtained with
a continuation path that starts from $\lambda=0$ and gradually grows
to $\lambda\approx 11.75$. At $\lambda=0$, the initial guess for
$u(x)$ is chosen to be a localized Gaussian function with power
$P=4$. For the following $q$ values, the solution of $u(x)$ of the
previous $\lambda$ value is used as the initial guess of the Newton
solve. Figure \ref{fig:t4}(b) plots for each $\lambda$ value the
number of Newton iterations as well as the averaged number of
preconditioned GMRES iterations for each Newton step. Notice that,
throughout the computation, the number of Newton iterations remains
rather small and this demonstrates the effectiveness of the Newton's
method and the continuation path. For each Newton step, the number of
the preconditioned GMRES iterations stays quite small and this is a
clear evidence of the effectiveness of the sparsifying
preconditioner. For our implementation, the setup cost of the
sparsifying preconditioner is about $4$ seconds and each application
of the sparsifying preconditioner takes about $0.02$ seconds. Figure
\ref{fig:t4}(c) and (d) give the field profile $u(x)$ at $\lambda=0$
and $\lambda=11.7498$, respectively. The last plot shows that the
support of the solution becomes quite extended as the eigenvalue
$\lambda$ approaches the bottom boundary of the lowest energy band of
the linear operator.

\begin{figure}[ht!]
  \begin{tabular}{cc}
    \includegraphics[height=1.5in]{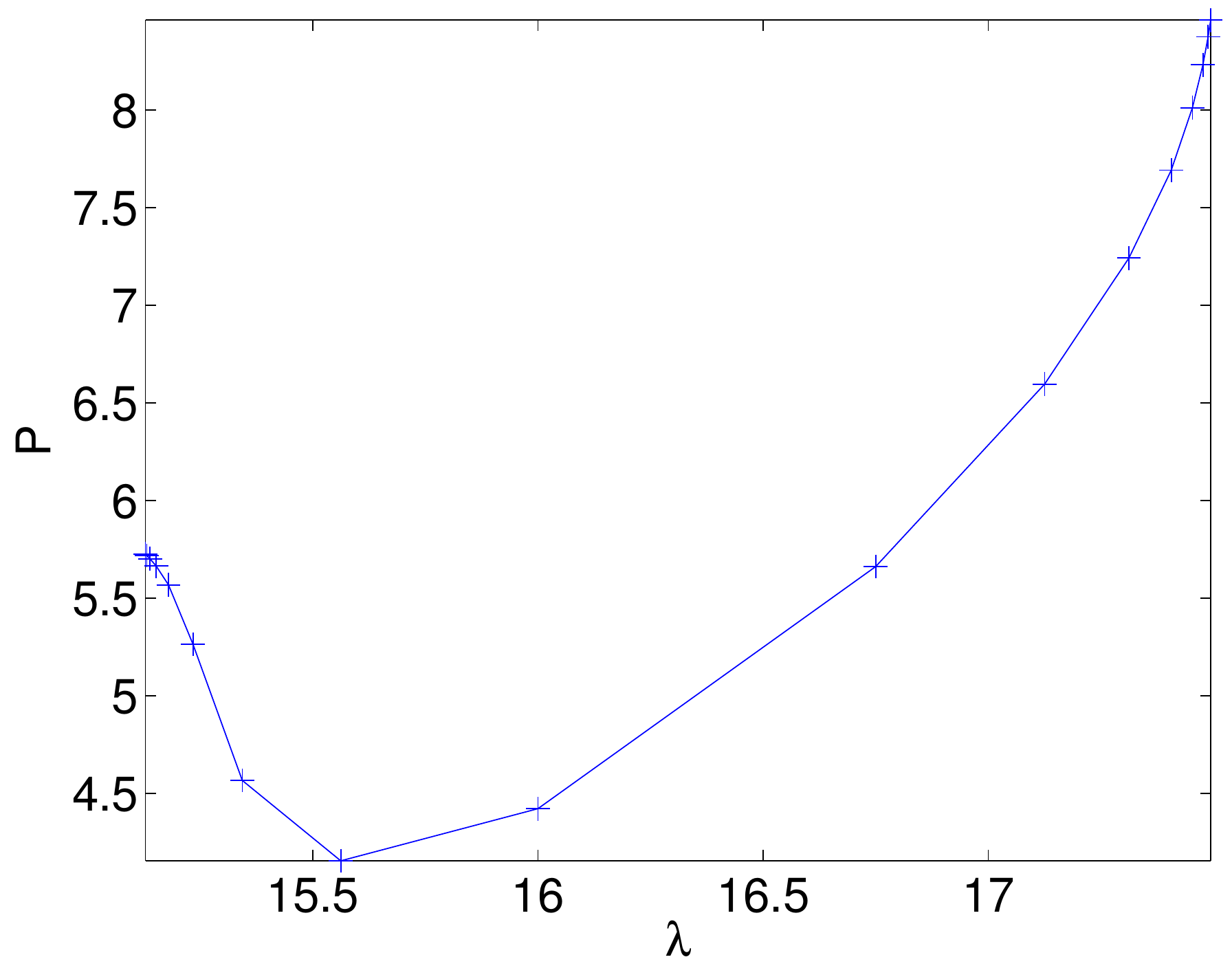} &    \includegraphics[height=1.5in]{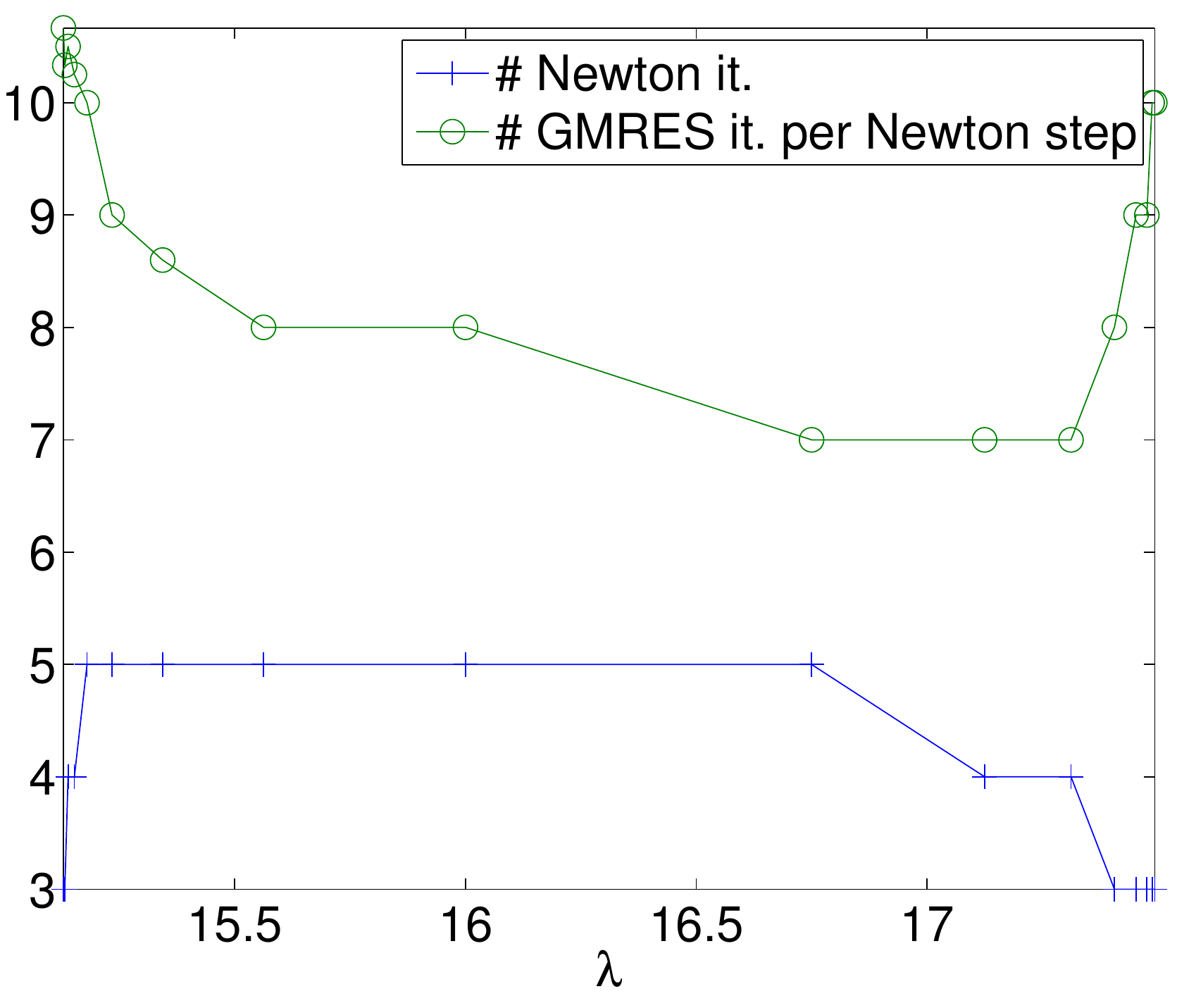} \\
    (a) & (b)\\
    \includegraphics[height=2in]{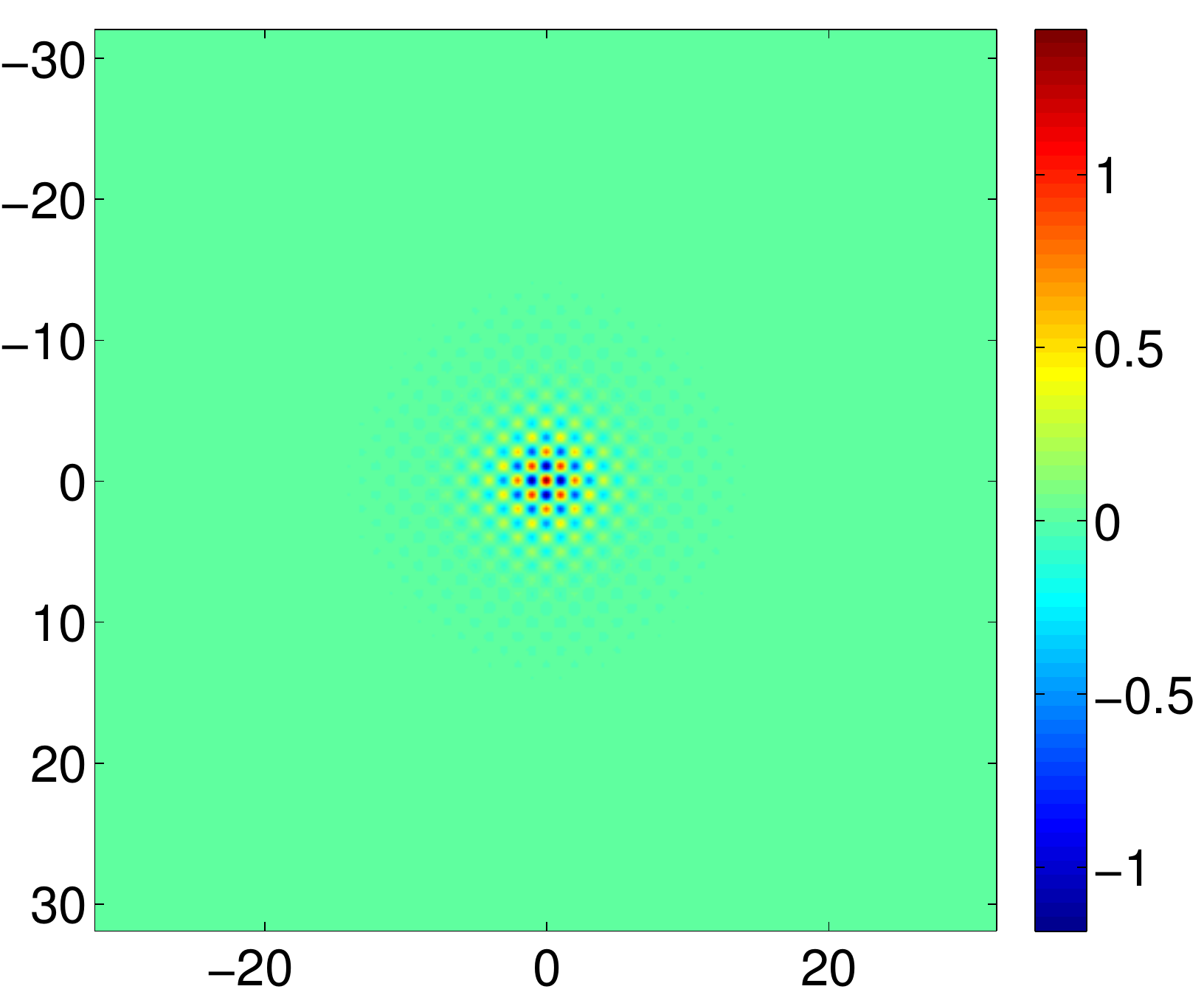} &     \includegraphics[height=2in]{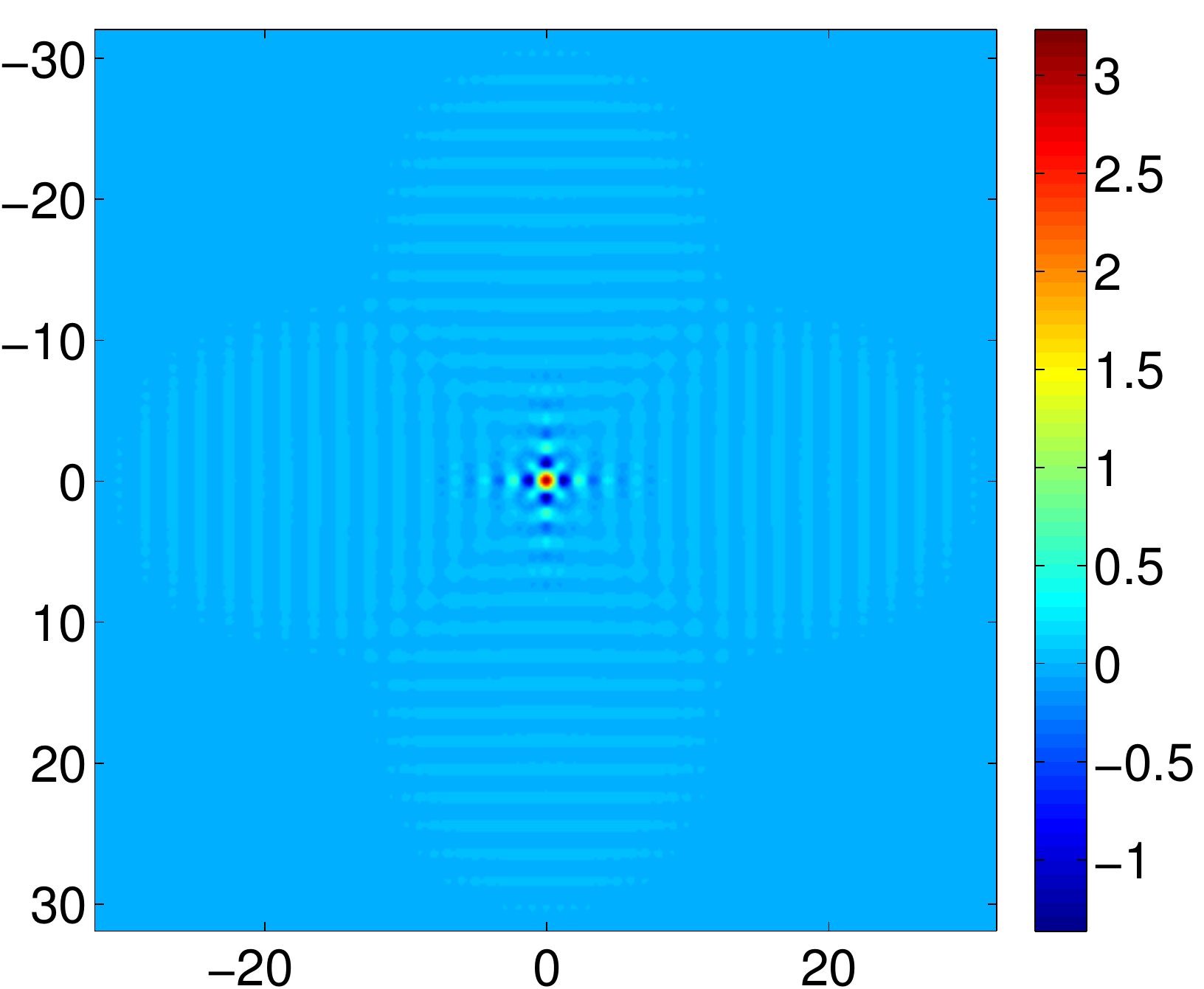}\\
    (c) & (d)
  \end{tabular}
  \caption{Kerr nonlinearity defocusing case $V_0=21.6$. (a) the
    $\lambda-P$ curve.  (b) The number of Newton iterations for each
    $\lambda$ value and the averaged number of preconditioned GMRES
    iteration for each Newton step.  (c) the field profile at
    $\lambda=15.1284$. (d) the field profile at $\lambda=17.4941$.  }
  \label{fig:t5}
\end{figure}

Figure \ref{fig:t5} shows the results of \eqref{eq:kerr} in the
defocusing case ($\sigma = -1$) with $V_0=21.6$. Since the support of
the defocusing gap soliton is significantly more extended, the
computation domain $[-32,32]^2$ is used here with $6$ points per unit
length. Figure \ref{fig:t5}(a) gives the $\lambda-P$ curve for
$\lambda\in (15.125, 17.5)$. This curve is obtained with two separated
continuation paths, one from $\lambda=16$ down to $\lambda=15.125$ and
the other from $\lambda=16$ to $\lambda=17.5$. In both paths, the
initial guess of $u(x)$ for $\lambda=16$ is a localized Gaussian
function with power $P=4$. Figure \ref{fig:t5}(b) plots for each
$\lambda$ the number of Newton iterations as well as the averaged
number of preconditioned GMRES iterations for each Newton step. The
curves show that our combination of the Netwon's method, the
continuation approach, and the sparsifying preconditioner results a
efficient way to compute the gap solitons in the defocusing case. Due
to the increase of the computational domain (as compared to the
focusing case), the setup cost of the sparsifying preconditioner here
is about $18$ seconds and each application of the sparsifying
preconditioner takes about $0.1$ seconds. Figure \ref{fig:t5}(c) and
(d) show the field profile $u(x)$ at $\lambda=15.1284$ and
$\lambda=17.4941$. In both plots, the essential support of $u(x)$
becomes extended since the eigenvalue $\lambda$ approaches the
boundary of the nearby energy bands of the linear operator.

\begin{figure}[ht!]
  \begin{tabular}{ccc}
    \includegraphics[height=1.5in]{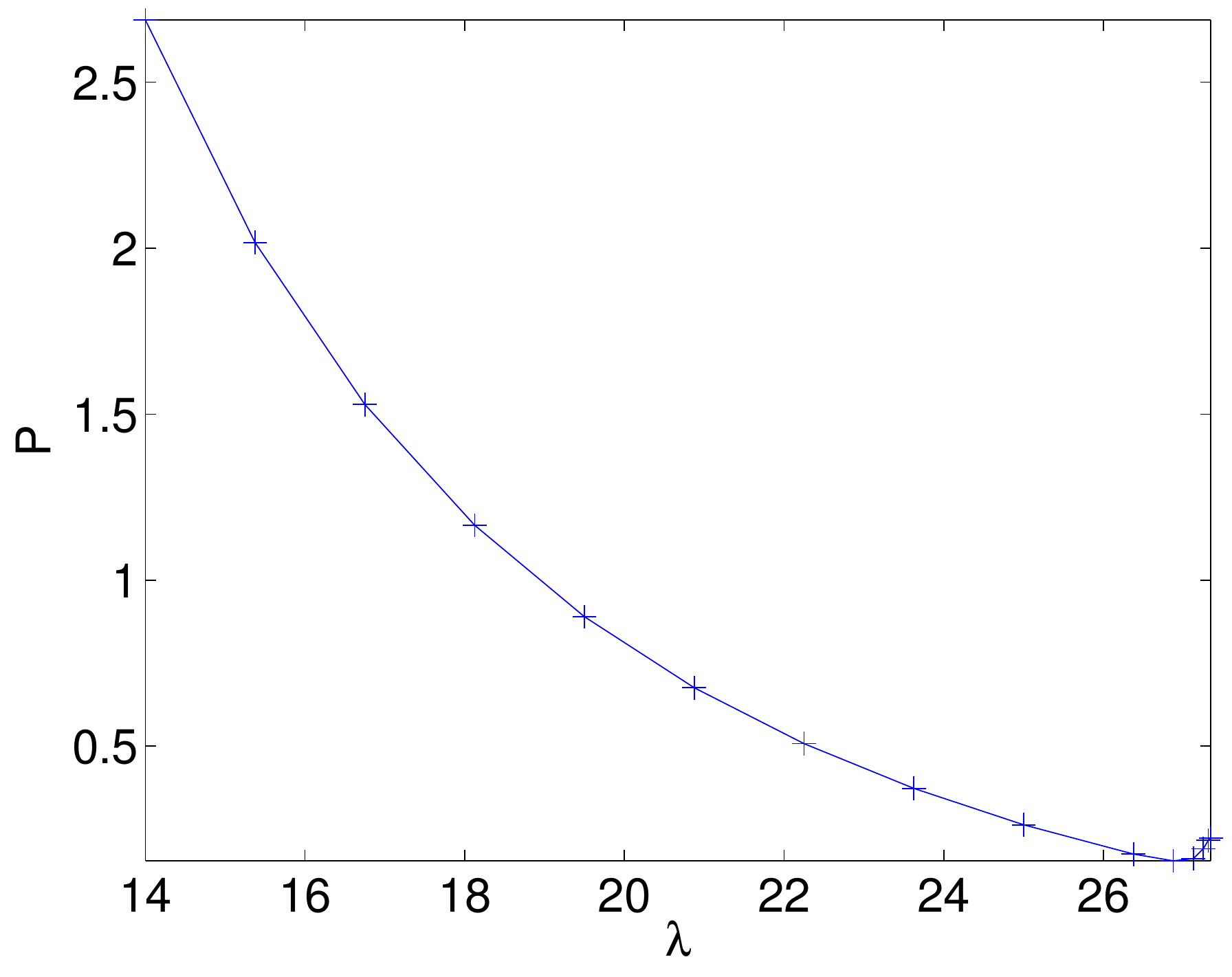}& \includegraphics[height=1.5in]{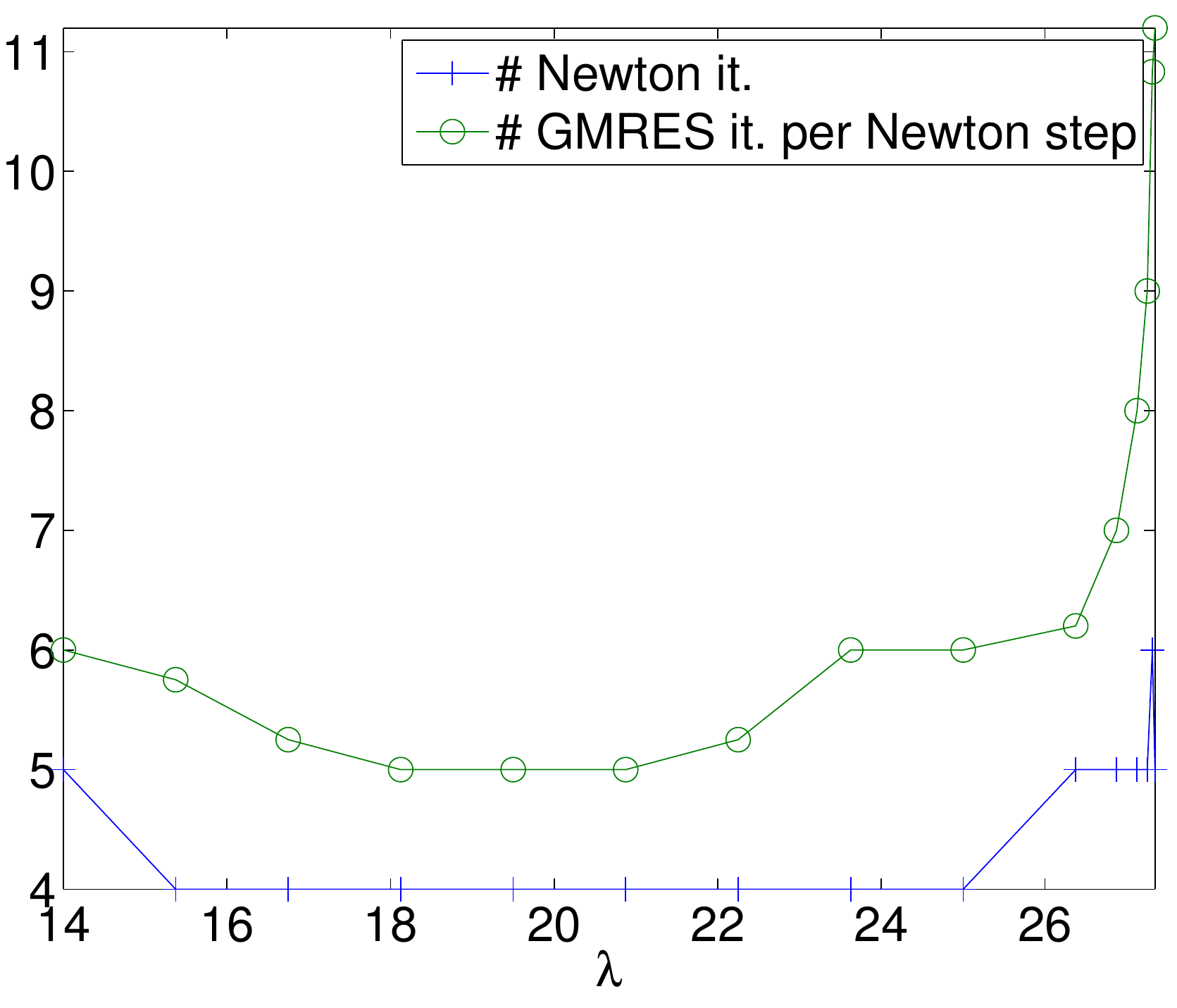}\\
    (a) & (b)\\
    \includegraphics[height=2in]{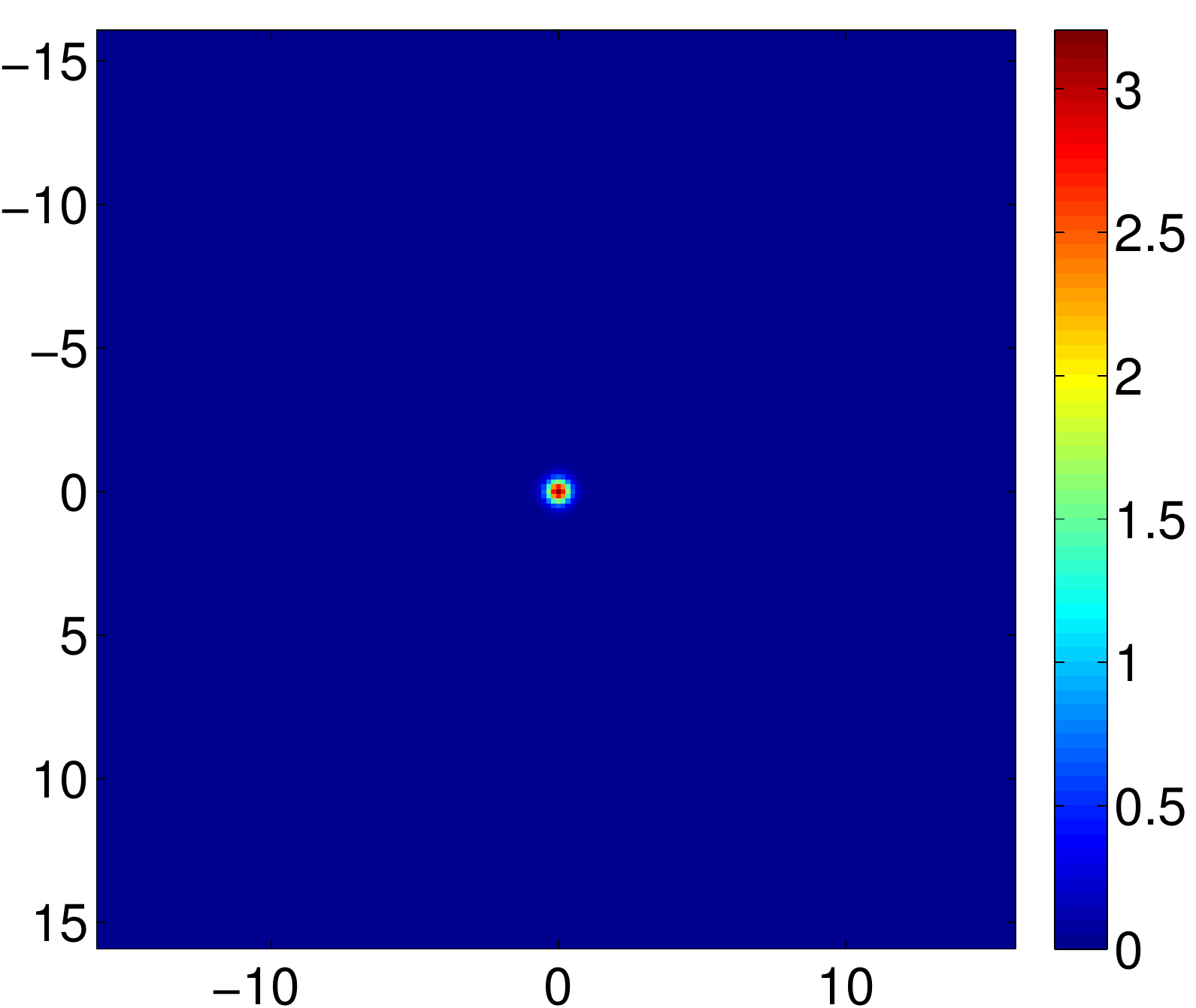}& \includegraphics[height=2in]{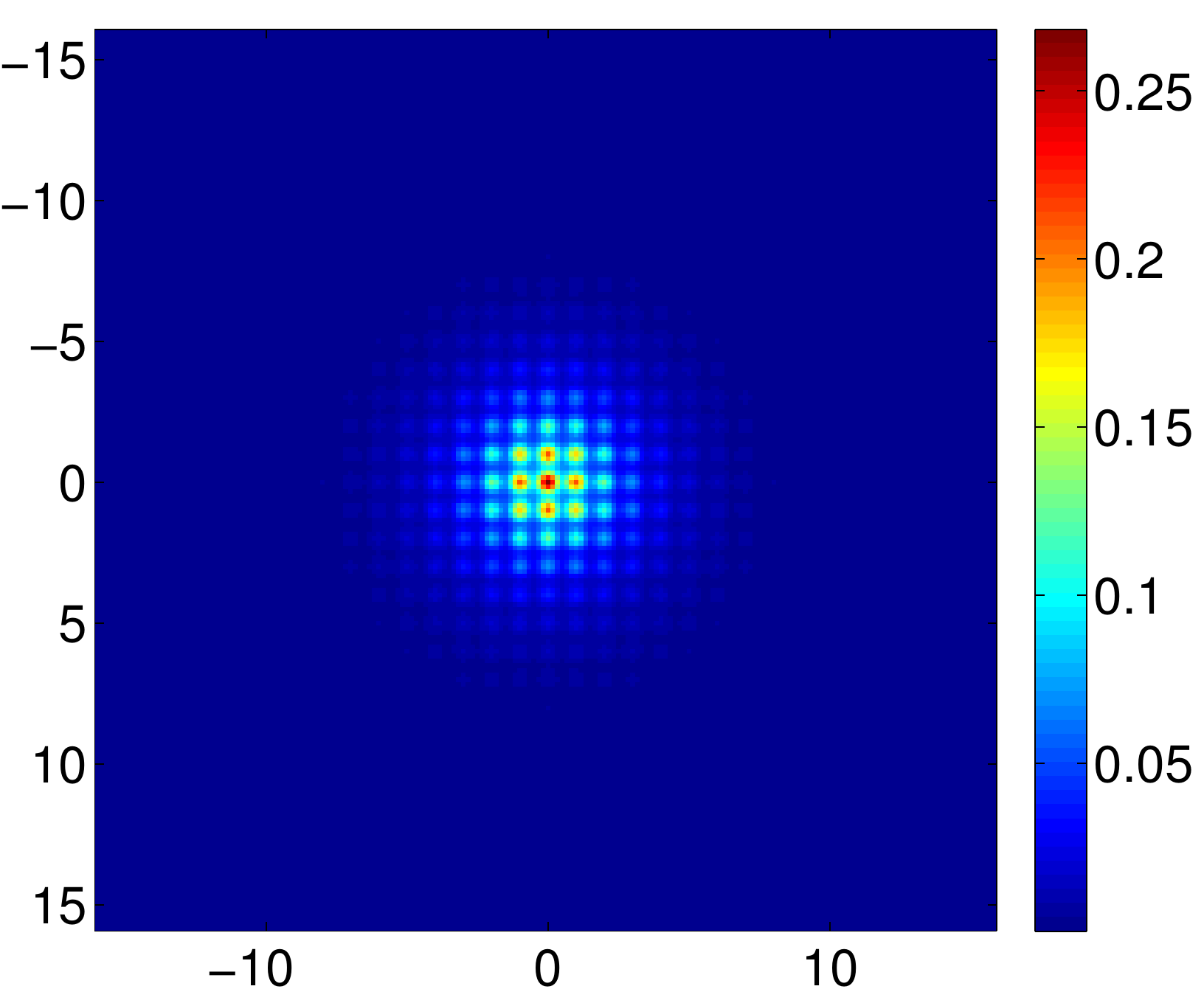}\\
    (c) & (d)
  \end{tabular}
  \caption{ Saturable nonlinearity focusing case $V_0=36.3$. (a) the
    $\lambda-P$ curve. (b) The number of Newton iterations for each
    $\lambda$ value and the averaged number of preconditioned GMRES
    iteration for each Newton step.  (c) the field profile at
    $\lambda=14$. (d) the field profile at $\lambda=27.3438$.  }
  \label{fig:t7}
\end{figure}

Figure \ref{fig:t7} summarizes the results of \eqref{eq:satu} in the
focusing case with $V_0=36.3$.  For this case, the computational
domain is set to be the periodic square $[-16,16]^2$ with $6$ points
per unit length. As before, Figure \ref{fig:t7}(a) is the $\lambda-P$
curve for $\lambda\in [14, 27.375]$. This curve is computed with a
single continuation path that starts from $\lambda=14$. At
$\lambda=14$, the initial guess for $u(x)$ is a localized Gaussian
profile with power $P=2$. Figure \ref{fig:t7}(b) gives for each
$\lambda$ value the number of Newton iterations and the averaged number
of preconditioned GMRES iterations for each Newton step. Figure
\ref{fig:t7}(c) and (d) plot the field profile at $\lambda=14$ and
$\lambda=27.3438$. Since the eigenvalue $\lambda=14$ is far from the
energy band, the field profile $u(x)$ in Figure \ref{fig:t7}(c) is
localized. On the contrary, $\lambda=27.3438$ approaches the lowest
energy band of the linear operator from bottom, the field profile
$u(x)$ is quite extended.

\begin{figure}[h!]
  \begin{tabular}{cc}
    \includegraphics[height=1.5in]{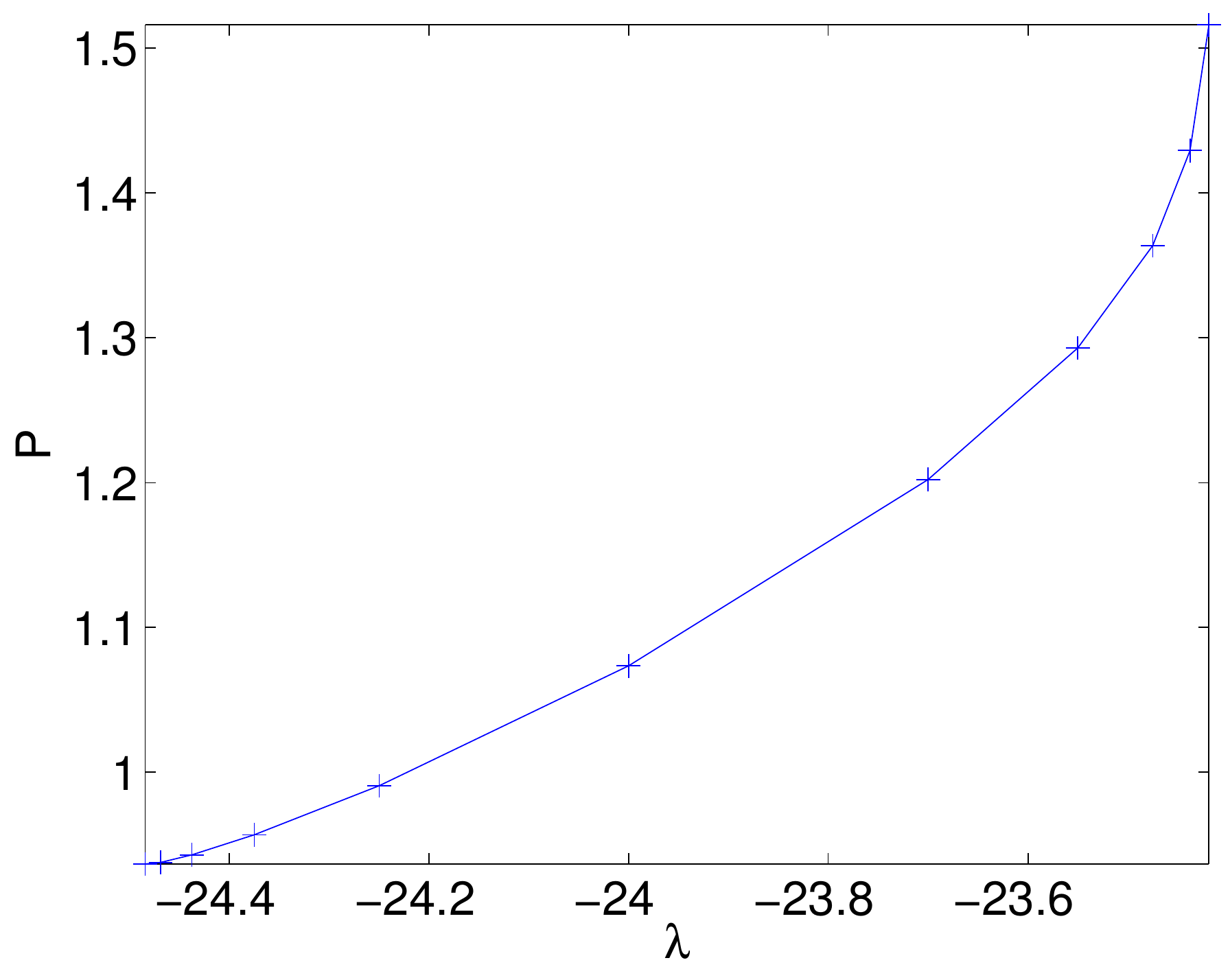}&\includegraphics[height=1.5in]{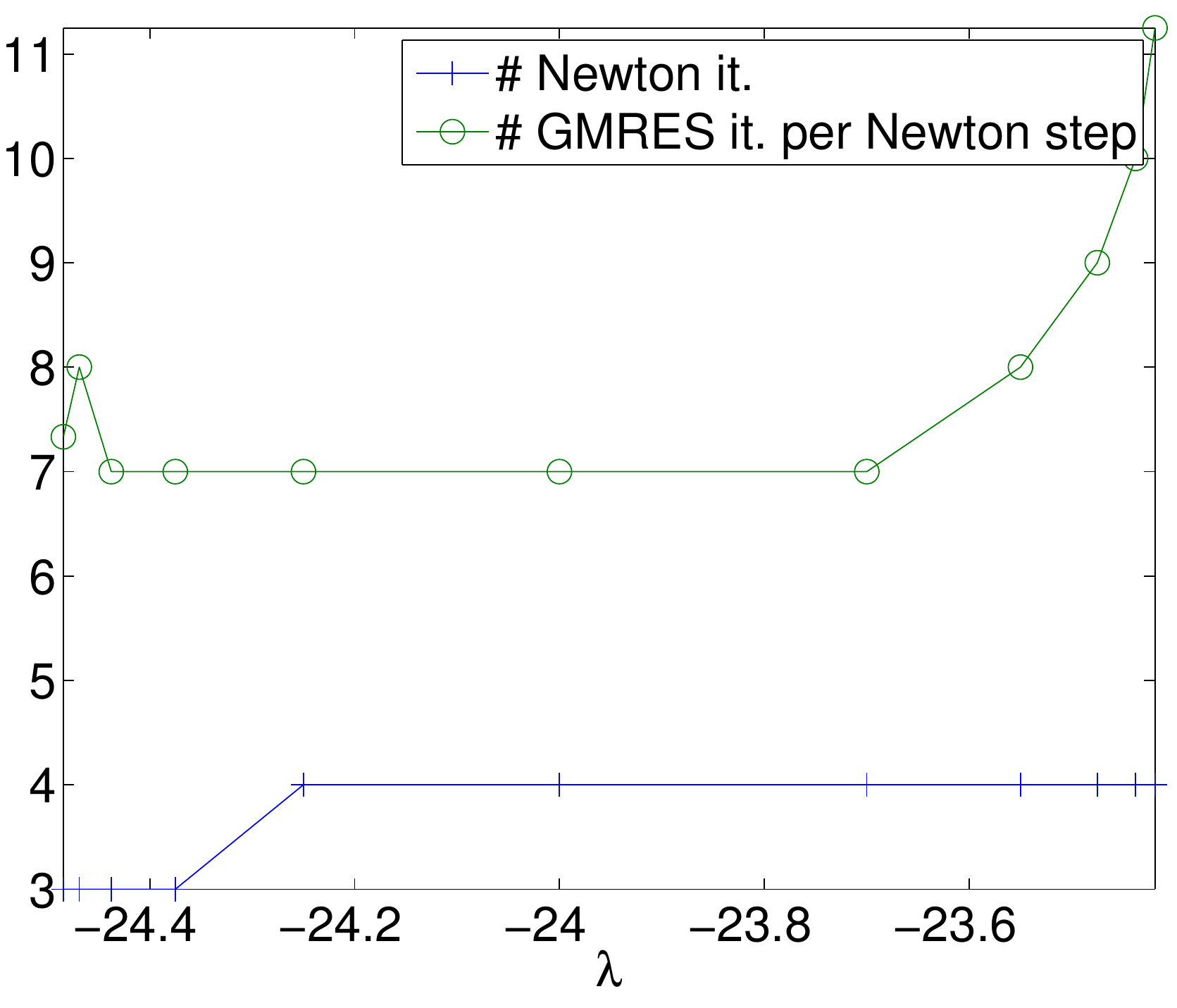}\\
    (a) & (b)\\
    \includegraphics[height=2in]{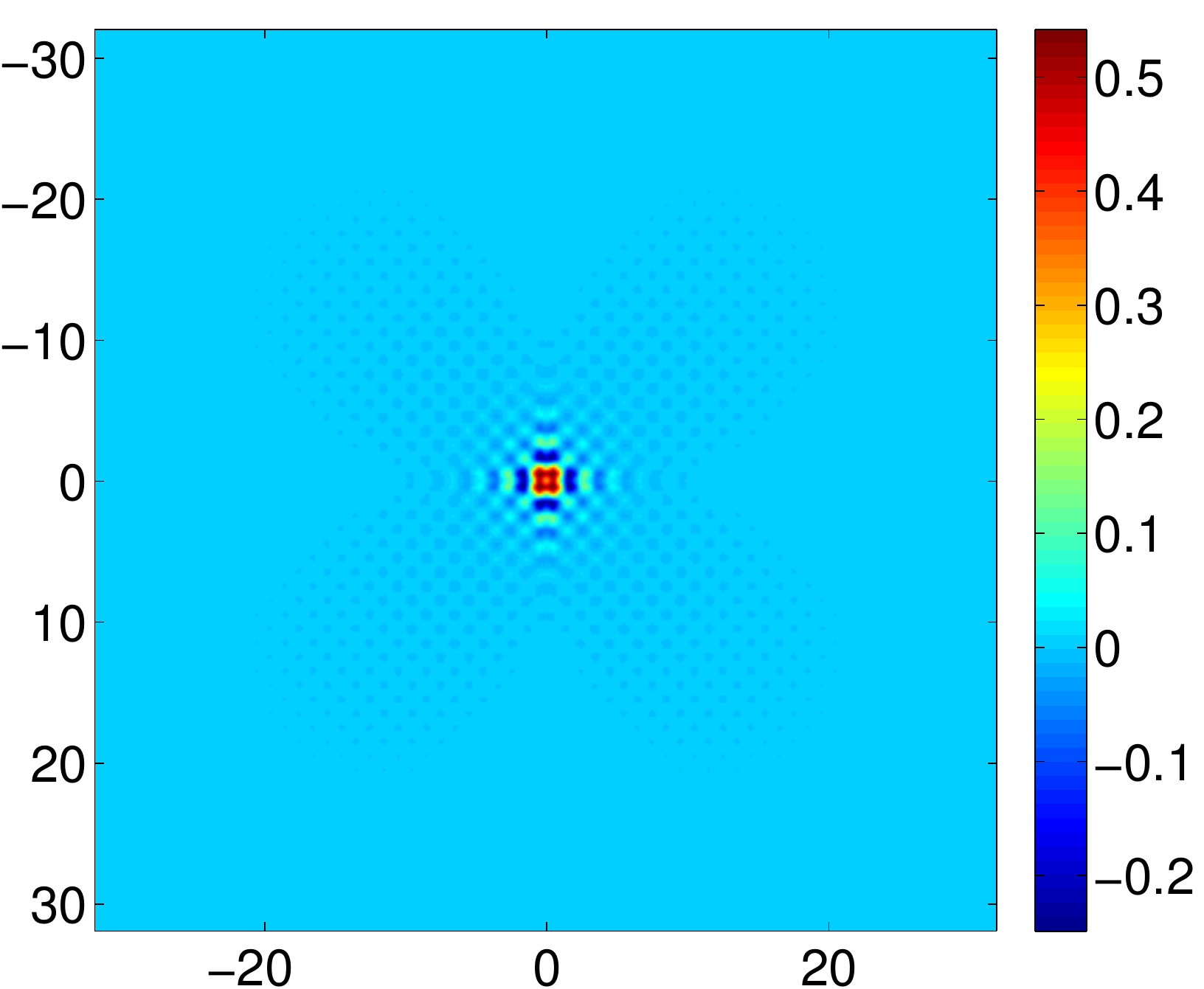}&    \includegraphics[height=2in]{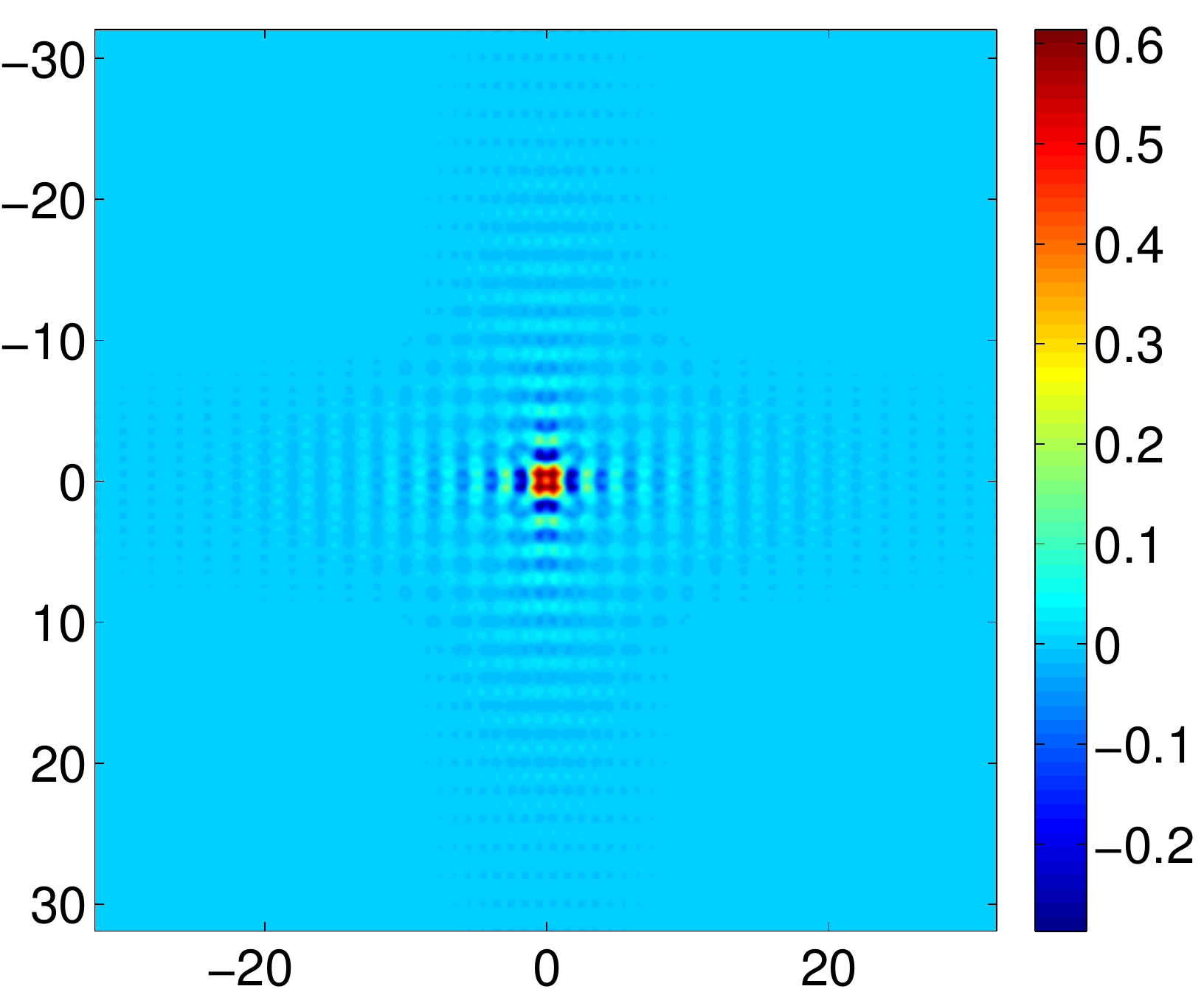}\\
    (c) & (d)
  \end{tabular}
  \caption{ Saturable nonlinearity defocusing case $V_0=-36.3$. (a)
    the $\lambda-P$ curve. (b) The number of Newton iterations for
    each $\lambda$ value and the averaged number of preconditioned
    GMRES iteration for each Newton step.  (c) the field profile at
    $\lambda=-24.4844$. (d) the field profile at $\lambda=-23.4187$.}
  \label{fig:t6}
\end{figure}

The results of \eqref{eq:satu} in the defocusing case with $V_0=-36.3$
are given in Figure \ref{fig:t6}. The computation domain is chosen to
be $[-32,32]^2$, due to the large support of the field profiles in the
defocusing case. As before, we discretize the domain with $6$ points
per unit length. Figure \ref{fig:t6}(a) plots the $\lambda-P$ curve
for $\lambda\in (-24.5, -23.4)$. This curve is obtained by following
two continuation paths, one from $\lambda=-24$ down to $\lambda=-24.5$
and the second from $\lambda=-24$ to $\lambda=-23.4$. For both paths,
the initial guess of $u(x)$ at $\lambda=-24$ is taken to be a
localized Gaussian profile with energy $P=0.4$. We notice that the
curve obtained is quite different from the one reported in
\cite{Efremidis:03}. Figure \ref{fig:t6}(b) plots the number of Newton
iterations and the averaged number of preconditioned GMRES iterations
for each Newton step, as a function of $\lambda$. Figure
\ref{fig:t6}(c) and (d) are the field profiles at $\lambda=-24.4844$
and $\lambda=-23.4187$, respectively. In both plots, the field profile
$u(x)$ is quite extended since the $\lambda$ value is close to the
boundary of the bottom and top energy bands of the linear operator,
respectively.

\section{Conclusion and Remarks}

In this work, we provide a robust and efficient method for soliton
calculations based on Newton iteration and sparsifying preconditioner
for the linearized problem. The Newton iteration based method is much
more robust than Pitviashvili type approaches.

As for possible extensions, while in this paper we solve the nonlinear
eigenvalue problem with prescribed eigenvalue $\lambda$ and
undetermined normalization, the sparsifying preconditioner can be also
applied to the problem with fixed normalization $\norm{u} = m$ and
undetermined $\lambda$. In the latter case, the Newton iteration becomes
\begin{equation}
  \begin{pmatrix}
    u^{n+1} \\
    \lambda^{n+1} 
  \end{pmatrix}
  = 
  \begin{pmatrix}
    u^n \\
    \lambda^n 
  \end{pmatrix}
  - 
  \begin{pmatrix} 
    -\Delta + L_{u^n} - \lambda^n & - u^n \\
    T_{u^n} & 0 
  \end{pmatrix}^{-1}
  \begin{pmatrix}
    - \Delta u^n + N(x, u^n) - \lambda^n u^n \\
    (\norm{u}_2^2 - m^2) / 2
  \end{pmatrix}.
\end{equation}
Here $L_u$ is the linearization of the potential terms at $u$: $L_u =
V + \frac{\delta N}{\delta u}$ and $T_u$ takes the inner product with
$u$: $T_u v = \average{u, v}$. In this case, the key step  is the solution to the linear system 
\begin{equation}
  \mc{L}_{(u, \lambda)} 
  \begin{pmatrix}
    v \\
    \mu 
  \end{pmatrix} 
  :=  \begin{pmatrix} 
    -\Delta + L_{u} - \lambda & - u \\
    T_{u} & 0 
  \end{pmatrix} 
  \begin{pmatrix}
    v \\
    \mu
  \end{pmatrix} 
  =   
  \begin{pmatrix}
    r \\
    \kappa
  \end{pmatrix}.
\end{equation}
We can apply sparsifying preconditioner to solve this equation, and
hence obtain an efficient algorithm.

Moreover, while in this work we focus on nonlinear Schr\"odinger
equations arising from optics application, the method can be applied
to other scenarios, for instance the nonlinear Maxwell equations.  In
addition, the sparsifying preconditioners also have potential
applications in the context of Kohn-Sham density functional theory
calculations. These directions will be explored in future works.

\bibliographystyle{amsxport}
\bibliography{nonSch}

\end{document}